\documentclass[11 pt]{amsart}
\usepackage{amssymb}
\usepackage{amsmath}
\usepackage{amsfonts}
\usepackage{graphicx}
\usepackage{amsthm}

\usepackage{enumerate}
\usepackage[mathscr]{eucal}
\usepackage{verbatim}

 \theoremstyle{plain}
\newtheorem{theorem}{Theorem}[section]
\newtheorem{lemma}[theorem]{Lemma}

\numberwithin{equation}{section}
\theoremstyle{plain}

\theoremstyle{remark}

\def\bbE{{\mathbb {E}}}
\def\bbR{{\mathbb {R}}}

\newcommand{\vertiii}[1]{{\left\vert\kern-0.25ex\left\vert\kern-0.25ex\left\vert #1 
    \right\vert\kern-0.25ex\right\vert\kern-0.25ex\right\vert}}

\begin{document}

\date{October, 2014}

\title
{Spherical means and pinned distance sets}
\author[]
{Daniel Oberlin and Richard Oberlin}

\address
{Daniel  Oberlin \\
Department of Mathematics \\ Florida State University \\
 Tallahassee, FL 32306}
\email{oberlin@math.fsu.edu}

\address
{Richard Oberlin \\
Department of Mathematics \\ Florida State University \\
 Tallahassee, FL 32306}
\email{roberlin@math.fsu.edu}

\subjclass{28E99}
\keywords{spherical averaging operator, Hausdorff dimension, pinned distance set}

\thanks{D.O. was supported in part by NSF Grant DMS-1160680
and R.O. was supported in part by NSF Grant DMS-1068523.}

\begin{abstract}
We use mixed norm estimates for the spherical averaging operator to obtain some results concerning pinned distance sets. 
\end{abstract}

\maketitle




\section{Introduction}

Let $E$ be a Borel subset of $\bbR^d$ and suppose $x\in\bbR^d$. Define 
\begin{equation*}
D_x (E)=\{|x-y|: y\in E\}.
\end{equation*}
Then $D_x (E)$ is called a {\it pinned distance set} of $E$. If the Hausdorff dimension 
$\dim E$ of $E$ is equal to $\beta$ and if $0\le\tau <1$, one may consider the possibility of estimates of the form 
\begin{equation}\label{1}
\dim  \{x\in\bbR^d: \dim D_x (E) < \tau \}\le \alpha 
\end{equation}
where $\alpha = \alpha (d,\beta ,\tau )$. For example, Theorem 8.3 in \cite{PS} shows that if $0\le \tau \le \min (1,\beta)$ then one may take $\alpha = d+\tau -\max(1,\beta )$.
One purpose of this note is to point out a relationship between estimates \eqref{1} and certain mixed norm estimates for the spherical averaging operator $S$ defined for 
nice functions $f$ on $\bbR^d$ and for $x\in\bbR^d$, $r>0$ by
\begin{equation*}
Sf(x,r)=\int_{S^{d-1}}f(x-r\sigma)\, d\sigma .
\end{equation*}
Here $d\sigma$ indicates integration with respect to Lebesgue measure on the unit sphere $S^{d-1}$ in $\bbR^d$. The mixed norm estimates we have in mind will be of the form 
\begin{equation}\label{2}
\Big(\int_{\bbR^d}\big(\int_{r_0}^{R_0}
|Sf(x,r)|^s\,dr\big)^{q/s}\,
d\lambda (x)\Big)^{1/q}\lesssim \vertiii {f } 
\end{equation}
where $\lambda$ is an $\alpha$-dimensional measure on $\bbR^d$ 
and where $\vertiii{\cdot}$ stands for either an $L^p$ or a Sobolev norm on $\bbR^d$. 

To describe the relationship between estimates \eqref{1} and \eqref{2} requires a little background. We begin by repeating some material from \cite{O}: 
for $\rho >0$, let $K_\rho$ be a kernel defined on $\bbR ^d$ by 
\begin{equation}\label{1.25}
K_\rho (x)=|x|^{-\rho}1_{B(0,r(d))}(x) 
\end{equation}
where $r=r(d)$ is a positive parameter.
Suppose that the Borel probability measure $\nu$ on $\bbR^d$ is a $\gamma$-dimensional measure 
in the sense that $\nu \big( B(x,\delta )\big)\lesssim \delta ^{\gamma}$ for all $x\in\bbR ^d$ 
and $\delta >0$. If $\rho <\gamma$ it follows that 
\begin{equation*}
\nu\ast K_\rho \in L^{\infty}(\bbR ^d ).
\end{equation*}
Also 
\begin{equation*}
\nu\ast K_\rho \in L^{1}(\bbR ^d )
\end{equation*}
so long as $\rho <d$. Thus, for $\epsilon>0$ and $1<p<\infty$, we have
\begin{equation}\label{ineq8}
\nu\ast K_\rho \in L^{p}(\bbR ^d ), \ \rho=\gamma +\frac{1}{p}(d-\gamma)-\epsilon
\end{equation}
by interpolation. 
The following lemma from \cite{O} is a weak converse of this observation.

\begin{lemma}\label{besov} If \eqref{ineq8} holds with $\epsilon =0$ and $1<p<\infty$, then 
$\nu$ is absolutely continuous with respect to Hausdorff measure of dimension
$\gamma -\epsilon '$ for any $\epsilon '>0$. Thus the support of $\nu$ has Hausdorff dimension 
at least $\gamma$.
\end{lemma} 
Returning to the relationship between \eqref{1} and \eqref{2}, suppose that $\nu$ is a Borel probability measure on $E\subset\bbR^d$. For each $x\in\bbR^d$ define the probability measure $\nu_x$ on $D_x (E)$ by 

\begin{equation*}
\int_{[0,\infty )} g \,d\nu_x =\int_E g(|x-e|)\, d\nu (e).
\end{equation*}

\noindent The proof of the next lemma will be given in \S 3.

\begin{lemma}\label{lemma2}
Suppose $\nu$ is as above, $\rho>d-1$, and $0<r_0 <R_0$. Then,
given $r(1)$ in \eqref{1.25},
it is possible to choose $r(d)$ (in \eqref{1.25}) so that 
we have the estimate
\begin{equation}\label{3}
\nu_x \ast K_{\rho+1-d}\,(r)\lesssim S(\nu\ast K_\rho )(x,r)\ (x\in\bbR^d ,\ r_0 \le r\le R_0 ) .
\end{equation}

\end{lemma}

\noindent Here, then, is a rough sketch (we have neglected $\epsilon$'s and various other details) of the argument which shows how estimates \eqref{2} for 
$\vertiii{\cdot}=\|\cdot \|_p$ 
can imply estimates \eqref{1}. (This argument is analogous to one in \cite{O} for orthogonal projections.) Suppose that $E\subset\bbR^d$ carries a $\beta$-dimensional probability measure $\nu$. Then
$\nu\ast K_\rho \in L^p (\bbR^d)$ for $\rho$ given by
\begin{equation*}
\rho=\beta+\frac{1}{p}(d-\beta ).
\end{equation*}
If \eqref{2} holds for 
$\vertiii{\cdot}=\|\cdot \|_p$, then 
\eqref{3} shows that $\nu_x \ast K_{\rho+1-d}\in L^s ([r_0 ,R_0 ])$ for $\lambda$ almost all
$x\in\bbR^d$. It will then follow from Lemma \ref{besov} that
for $\tau$ given by 
\begin{equation*}
\rho+1-d =\tau +\frac{1}{s}(1-\tau )
\end{equation*}
the support of $\nu_x$ has dimension at least $\tau$ for $\lambda$ almost all
$x$. Since $\nu_x$ is supported on 
$D_x (E)$ and since $\lambda$ is an arbitrary $\alpha$-dimensional measure, \eqref{1} follows.

The remainder of this note is organized as follows: \S 2 contains the statements of our results and \S 3 contains their proofs.

\section{Results}

Estimates such as \eqref{2} are examples of what Wolff \cite{W} called \lq\lq $L^p \rightarrow L^q$ inequalities for the wave equation relative to fractal measures" - see also 
\cite{M}, \cite{O2}, \cite{O3}. The following example, an estimate for a \lq\lq fractal spherical maximal  function", is an easy consequence of Theorem $4_S$ in \cite{O3}. 

\begin{theorem}\label{est1}
Suppose that $\lambda$ is a nonnegative compactly-supported Borel measure on $\bbR^d$  and suppose that for some $\alpha\in(0,(d-1)/2)$ the measure $\lambda$ satisfies an estimate
\begin{equation}\label{2.05}
\lambda \big(B(x,\delta )\big)\lesssim \delta^\alpha
\end{equation}
for $x\in\bbR^d$, $\delta>0$. Then, for $\epsilon >0$,
$q<2$, and $0<r_0 <R_0 <\infty$, there is the estimate 
\begin{equation}\label{2.1}
\Big(\int_{\bbR^d}\big(
\sup_{r_0 <r<R_0}
|Sf(x,r)|
\big)^{q}\,
d\lambda (x)\Big)^{1/q}\lesssim \|f\|_{W^{2,(1-\alpha )/2-\epsilon}}\ .
\end{equation}
The constant implicit in \eqref{2.1} depends on the size of the support of $\lambda$, the constant implicit in \eqref{2.05}, $\epsilon$, $q$, $r_0$, and $R_0$. 
\end{theorem}
The next estimate, of form \eqref{1}, is a consequence of Theorem \ref{est1}:
\begin{theorem}\label{est2}
Suppose that the compact set $E\subset\bbR^d$ satisfies $\dim E =\beta $ for some $\beta >0$. Suppose that $0<\tau <1$ and that $2\tau +(d-1)/2<\beta <2\tau +d-1$. Then 
\begin{equation*}
\dim  \{x\in\bbR^d: \dim D_x (E)< \tau \}\le 2\tau -\beta+d-1 .
\end{equation*}
\end{theorem}
%

The proof of Theorem $4_S$ in \cite{O3} depends on Fourier analysis. Theorem \ref{est3} below is proved by modifying a measure-theoretic argument previously used in \cite{O} and \cite{O2}.

\begin{theorem}\label{est3} 
Suppose $\lambda$ is a compactly supported probability measure on $\bbR^d$ and,
for any fixed $0<r_0 <R_0 <\infty$, consider the mixed norm estimate
\begin{equation}\label{2.5}
\Big(\int_{\bbR^2}\big(\int_{r_0}^{R_0} |Sf (x,r)|^s \,dr\big)^{q/s}\,d\lambda (x)\Big)^{1/q}
\lesssim \|f\|_p .
\end{equation}

\noindent(a) Suppose $d=2$ and $\alpha >1/2$. Suppose that $\lambda$ satisfies the Frostman condition 
\begin{equation}\label{2.3}
\int_{\bbR^2}\int_{\bbR^2}\frac{1}{|x_1 -x_2 |^\alpha}\,d\lambda (x_2 )\,d\lambda (x_1 )\le C<\infty .
\end{equation}
Then \eqref{2.5} holds when 
\begin{equation*}
\Big(\frac{1}{p},\frac{1}{q},\frac{1}{s}\Big)=t\Big(\frac{1}{2},\frac{1}{2\alpha},\frac{1}{4}\Big)
+(1-t)\Big(1,0,1\Big),\  0\le t<1 .
\end{equation*}

\noindent(b) Suppose $d=2$ and $0<\alpha <\alpha ' < 1/2$. Suppose that $\lambda$ satisfies the condition 
\begin{equation}\label{2.6}
\lambda \big( B(x,\delta )\big)\le C \delta^{\alpha '} \  (x\in\bbR^2, \delta>0). 
\end{equation}
Then \eqref{2.5} holds where 
\begin{equation}\label{2.65}
\Big(\frac{1}{p},\frac{1}{q},\frac{1}{s}\Big)=t\Big(\frac{1}{1+2\alpha },\frac{1}{1+2\alpha },\frac{1-\alpha }{1+2\alpha }\Big)
+(1-t)\Big(1,0,1\Big),\  0\le t<1 .
\end{equation}

\noindent(c) Suppose $d>2$ and $0<\alpha <\alpha ' < 1$. Suppose that $\lambda$ satisfies the condition 
\begin{equation}\label{2.66}
\lambda \big( B(x,\delta )\big)\le C \delta^{\alpha '} \  (x\in\bbR^d, \delta>0). 
\end{equation}
Then \eqref{2.5} holds where 
\begin{equation}\label{2.67}
\Big(\frac{1}{p},\frac{1}{q},\frac{1}{s}\Big)=t\Big(\frac{1}{1+\alpha },\frac{1}{1+\alpha },\frac{1-\alpha }{1+\alpha }\Big)
+(1-t)\Big(1,0,1\Big),\  0\le t<1 .
\end{equation}
\end{theorem}

\noindent The implied constant in \eqref{2.5} depends on $r_0$, $R_0$, the size of the support of $\lambda$, $C$  
from \eqref{2.3}, \eqref{2.6}, or \eqref{2.66}, and $t$. There is an analogue of (c) of Theorem \ref{est3} for $\alpha >1$, but it is essentially subsumed by Theorem \ref{est2}.
A consequence of Theorem \ref{est3} is the following analogue of Theorem 1.2 in \cite{O}:

\begin{theorem}\label{est4}
Suppose $E\subset \bbR^d$ has Hausdorff dimension $\beta$. 

\noindent(a) If $d=2$ and $\beta >1/2$, then 
\begin{equation*}
\dim\{x\in\bbR^2:\dim D_x (E) < (2\beta -1)/3\}=0 .
\end{equation*}

\noindent(b) If $d>2$ and $\beta >d-2$, then 
\begin{equation*}
\dim\{x\in\bbR^2:\dim D_x (E) < (\beta +2-d)/2\}=0 .
\end{equation*}

\end{theorem}

\section{Proofs}

{\it Proof of Lemma \ref{lemma2}:} Suppose $\delta \in (0,r_0 /2 )$. Write $\sigma_{x,r}$ for normalized Lebesgue measure on 
$\{y\in\bbR^d :|x-y|=r\}$. Then if $ r_0 \le r\le R_0$ we have
\begin{multline*}
\nu_x ([r-\delta ,r+\delta ])=\int_E 1_{[r-\delta ,r+\delta ]}(|x-e |)\,d\nu (e) \approx \\
\langle \delta^{1-d}1_{B(0,\delta )}\ast \sigma_{x,r},\nu \rangle = 
 \langle \sigma_{x,r},\delta^{1-d}1_{B(0,\delta )}\ast\nu\rangle =
S(\delta^{1-d}1_{B(0,\delta )}\ast\nu)(x,r),
\end{multline*}
where the implied constant depends on $R_0$. If $r(1)\le r_0 /4$ and 
\begin{equation*}
K_{\rho+1-d}(x)=|x|^{-(\rho+1-d)}1_{B(0,r(1))}(x),\ x\in\bbR
\end{equation*}
then 
\begin{multline*}
\nu_x \ast K_{\rho+1-d}(r) \approx \sum_{2^{-j}\le r(1)}2^{(\rho +1-d)j}\nu_x ([r-2^{-j},r+2^{-j}])=\\
\sum_{2^{-j}\le r(1)}2^{(\rho +1-d)j}S\big((2^{-j})^{-(d-1)}1_{B(0,2^{-j})}\ast \nu \big) (x,r) \lesssim 
S(K_\rho \ast \nu )(x,r) ,
\end{multline*}
so long as $r(d)\gtrsim r(1)$.

$\ \ \ \ \ \ \ \ \ \ \ \ \ \ \ \ \ \ \ \ \ \ \ \ \ \ \ \ \ \ \ \ \ \ \ \ \ \ \ \ \ \ \ \ \ \ \ \ \ \ \ \ \ \ \ \ \ \ \ \ \ \ \ \ \ \ \ \ \ \ \ \ \ \ \ \ \ \ \ \ \ \ \ \ \ \ \ \ \ \ \square$

\ \ \ 

{\it Proof of Theorem \ref{est1}:} As mentioned in \S 2, the proof of Theorem \ref{est1} is consequence of the following result from \cite{O3}:

\begin{theorem}\label{prev}
Suppose $\mu$ is a nonnegative Borel measure on a compact subset of $\bbR^d \times (0,\infty )$ and suppose that, for some $\alpha\in \big( 0,(d-1)/2\big)$, $\mu$ satisfies the estimate 
\begin{equation}\label{3.1}
\mu\big(\{ (x,r)\in \bbR^d \times (0,\infty):|x-x' |+|r-r' |<\delta\}\big)\lesssim \delta^\alpha
\end{equation}
for all $(x' ,r' )\in \bbR^d \times (0,\infty)$ and $\delta >0$. Then, for $\epsilon >0$, 
\begin{equation}\label{3.2} 
\|Sf\|_{L^{2,\infty}_\mu} \lesssim \|f\|_{W^{2,(1-\alpha )/2+\epsilon}}.
\end{equation}
\end{theorem}
\noindent Since $\mu$ is compactly supported, $\text{supp} (\mu )\subset B(0,M)\times [r_0 ,R_0]$
for some $M>0$ and $0<r_0 <R_0 <\infty$. 
 The proof of Theorem \ref{prev} shows that the constant implicit in \eqref{3.2} is, for fixed 
 $M$, $r_0$, $R_0$, 
 $\alpha$, and $\epsilon$, bounded by a function of the constant implicit in \eqref{3.1}. 
Suppose that $\lambda$ is as in the statement of Theorem \ref{est1} and  $x\mapsto r(x)$ is any Borel function from the compact support of $\lambda$ into $[r_0 ,R_0 ]$. 
If the measure $\mu$ on $\bbR^d \times (0,\infty)$ is defined by 
\begin{equation*}
\int_{\bbR^d \times (0,\infty)}f\,d\mu=\int_{\bbR^d}f\big( x,r(x)\big)\, d\lambda (x)
\end{equation*}
then the hypothesis \eqref{2.05} on $\lambda$ ensures that the measures $\mu$ 
satisfy \eqref{3.1} uniformly in the choice of $r(x)$. Thus 
\begin{equation*}
\|Sf \big(\cdot ,r(\cdot )\big)\|_{L^{2,\infty}_\lambda }  \lesssim \|f\|_{W^{2,(1-\alpha )/2+\epsilon}}
\end{equation*}
with implicit constant independent of $r(x)$. That is enough to establish \eqref{2.1}.
$\ \ \ \ \ \ \ \ \ \ \ \ \ \ \ \ \ \ \ \ \ \ \ \ \ \ \ \ \ \ \ \ \ \ \ \ \ \ \ \ \ \ \ \ \ \ \ \ \ \ \ \ \ \ \ \ \ \ \ \ \ \ \ \ \ \ \ \ \ \ \ \ \ \ \ \ \ \ \ \ \ \ \ \ \ \ \ \ \ \ \square$

\ \ \ \ 

{\it Proof of Theorem \ref{est2}:} Define $\alpha$ by $\alpha =2\tau -\beta +d-1$ and let $\lambda$ be any measure satisfying \eqref{2.05}. It will be enough to show that if $\tilde F\subset\bbR^d$ is any compact set with $\lambda (\tilde F)>0$ then there is a compact subset $F$ of $\tilde F$ 
with $\lambda (F)>0$
such that 
$\dim D_x (E)\ge\tau$ for $\lambda$-a.a.$\,x\in F$. 

Fix a small $\epsilon >0$ and set $\beta '=\beta -\epsilon$. 
%
%
Since $\dim E =\beta$ there is a probability measure $\nu$ on $E$ with 
\begin{equation}\label{3.35}
\int_{\bbR^d}\frac{|\hat\nu (\xi )|^2}{|\xi |^{d-\beta '}}\, d\xi <\infty .
\end{equation}

Now fix $\tilde F$ with $\lambda (\tilde F) >0$ as above and choose $0<r_0 <R_0$ (in \eqref 2), $r(1)$ (in the definition of the one-dimensional kernel $K_\rho$), and $F\subset\tilde F$ with $\lambda (F)>0$ such that $\nu_x ([r_0 +r(1),R_0 -r(1)])>0$ for $\lambda$-a.a.$\,x\in F$. 
Choose $r(d)$ in the definition of the $d$-dimensional kernel $K_\rho$ so that \eqref{3} holds.
Since the $d$-dimensional $K_\rho$ satisfies $|\hat K_\rho (\xi )|\lesssim |\xi |^{-d+\rho}$ it follows from \eqref{3.35} that 
$K_\rho \ast \nu \in W^{2,(1-\alpha)/2 -\epsilon /2}$ if 
$$
\rho = (\alpha +\beta +d-1)/2.
$$
Then Theorem \ref{est1} and Lemma \ref{lemma2} show that for any $s<\infty$ we have, for the one-dimensional kernel $K_{(\alpha+\beta+1-d)/2}$, that 
$$
\nu_x \ast K_{(\alpha+\beta+1-d)/2}\in L^s (r_0 ,R_0 )
$$ 
for $\lambda$-a.a.$\,x\in F$. It follows from Lemma \ref{besov} (with $d=1$) and the fact that $\nu_x$ is supported on $D_x (E)$ that 
\begin{equation*}
\dim D_x (E)\ge \frac{(\alpha +\beta +1-d)/2-1/s}{1-1/s}=\frac{\tau-1/s}{1-1/s}
\end{equation*}
for every $s<\infty$ for $\lambda$-a.a.$\,x\in F$. Letting $s\rightarrow\infty$, we have $\dim D_x (E)\ge\tau$ for $\lambda$-a.a.$\,x\in F$ as desired.
$\ \ \ \ \ \ \ \ \ \ \ \ \ \ \ \ \ \ \ \ \ \ \ \ \ \ \ \ \ \ \ \ \ \ \ \ \ \ \ \ \ \ \ \ \ \ \ \ \ \ \ \ \ \ \ \ \ \ \ \ \ \ \ \ \ \ \ \ \ \ \ \ \ \ \ \ \ \ \ \ \ \ \ \ \ \ \ \ \ \ \square$

\ \ \ \ 

{\it Proof of Theorem \ref{est3}:} It is clear that \eqref{2.5} holds when $p=1$, $q=\infty$, and $s=1$.
Thus it is enough to prove a restricted weak type version of the result that would correspond to \eqref{2.5} with $t=1$. 
So suppose $E\subset\bbR^d$ and $S1_E \, (x,r)\geq \mu$ for
$$
(x,r)\in F=\{(x,r): x\in A\subset\bbR^d, \, r\in T_x\subset [r_0 ,R_0 ]\}
$$
where the one-dimensional measure $|T_x|$ satisfies $B\le |T_x|\le 2B$ for some positive $B$.
We will show that if $1/p,\,1/q/,\, 1/s$ correspond to $t=1$ then the inequality
\begin{equation}\label{4}
\mu\,\lambda (A)^{1/q}B^{1/s}\lesssim |E|^{1/p}
\end{equation}
holds. The strategy is to estimate $|E|$ from below. If $x\in \bbR^d$ let 
\begin{equation*}
E_x =\{x+r\omega\in E: r\in T_x,\, \omega\in S^{d-1} \}.
\end{equation*}
Then, for $x_1 ,\dots ,x_N \in A$ we have 
\begin{equation}\label{4.75}
|E|\ge \sum_{j=1}^N |E_{x_j}|-\sum_{1\le j<k\le N}|E_{x_j}\cap E_{x_k}|.
\end{equation}
It follows from our assumptions that if $x\in A$ then
\begin{equation}\label{4.5}
|E_x |\gtrsim B\mu. 
\end{equation}
To make use of this by way of the estimate \eqref{4.75} we will need an upper 
bound for $|E_{x_1}\cap E_{x_2}|$. Thus we will begin the proof of (a) by establishing the following (two-dimensional) estimate 
\begin{equation}\label{6}
|E_{x_1}\cap E_{x_2}|\lesssim \frac{B^{3/2}}{|x_1 -x_2 |^{1/2}}
\end{equation}
(where the implied constant depends on $r_0$ and the size of the support of $\lambda$).
{\it Proof of \eqref{6}}:
Assume without loss of generality that the $x_i$ lie on the $x$-axis and abuse notation by writing $x_i =(x_i ,0)$. Consider the transformation of $\bbR^2$ defined by 
\begin{equation*}
(y_1 ,y_2)\mapsto \big((x_1 -y_1 )^2 +y^2_2, (x_2 -y_1 )^2 +y_2^2 \big)\dot= (r_1^2 ,r_2^2 ).
\end{equation*}
We will establish \eqref{6} by showing that 
\begin{equation*}
\int_{T_{x_2}}\int_{T_{x_1}}\Big|\frac{\partial (y_1 ,y_2 )}{\partial (r_1 ,r_2 )}\Big|\, dr_1 \,dr_2
\lesssim \frac{B^{3/2}}{|x_1 -x_2 |^{1/2}}\ .
\end{equation*}
Since $r_0 \le r_1 ,r_2 \le R_0$ it is enough to show that 
\begin{equation}\label{8}
\int_{T_{x_2}}\int_{T_{x_1}}\Big|\frac{\partial (y_1 ,y_2 )}{\partial (r_1^2 ,r_2^2 )}\Big|\, dr_1 \,dr_2
\lesssim \frac{B^{3/2}}{|x_1 -x_2 |^{1/2}}
\end{equation}
holds.
A computation shows that 
\begin{equation*}
\Big|\frac{\partial (y_1 ,y_2 )}{\partial (r_1^2 ,r_2^2 )}\Big|=
\frac{1}{|y_2 |\, |x_1 -x_2 |}.
\end{equation*}
Now $ {|y_2 | \, |x_1 -x_2 |}$ is just twice the area of the triangle with vertices $(x_1 ,0), (x_2 ,0), (y_1 ,y_2 )$. Write $\theta$ for the angle between the $x$-axis and the line through 
$(x_1 ,0)$ and  $(y_1 ,y_2 )$ and $\Delta$ for $|x_1 -x_2 |$. From 
\begin{equation*}
r_1 \cos \theta =\frac{r_1^2 +\Delta ^2 -r_2^2}{2\Delta}
\end{equation*}
it follows that 
\begin{equation*}
2{|y_2 | \, |x_1 -x_2 |}=2\Delta r_1 \sin\theta=\sqrt{\big(r_2^2 -(r_1 -\Delta )^2\big)
\big(r_2^2 -(r_1 +\Delta )^2\big)} \ .
\end{equation*}
Since $r_1 ,r_2 \ge r_0$ we can estimate 
\begin{equation*}
\sqrt{\big(r_2^2 -|r_1 -\Delta |^2\big)\big(r_2^2 -|r_1 +\Delta |^2\big)} \gtrsim 
\sqrt{\big|\big(r_2 -|r_1 -\Delta |\big)
\big(r_2 -|r_1 +\Delta |\big)\big|}\  . 
\end{equation*}
Therefore \eqref{8} will follow from the estimate 
\begin{equation*}
\int_{T_{x_1}}\int_{T_{x_2}}\frac {1}{\sqrt{\big|\big(r_2 -|r_1 -\Delta |\big)\big(r_2 -|r_1 +\Delta |\big)\big|}}\,dr_1 \, dr_2 \lesssim \frac{B^{3/2}}{\Delta^{1/2}}\ .
\end{equation*}
Since $B\le |T(x_i )|\le 2B$, a scaling argument shows that this will follow from 
\begin{equation}\label{9}
\int_{T_1}\int_{T_2}\frac {1}{\sqrt{\big|\big(r_2 -|r_1 -1 |\big)\big(r_2 -|r_1 +1 |\big)\big|}}\,dr_1 \, dr_2 \lesssim {B^{3/2}}.
\end{equation}
where 
$B\le |T_i |\le 2B$ and
$T_i \subset (\eta ,\infty )$.
Here $\eta >0$ depends on $r_0$ and an upper bound for $\Delta =|x_1 -x_2 |$ and so on the size of the support of $\lambda$. 
Also the implied constant in \eqref{9} depends on $\eta$. To see \eqref{9}, note that, since $r_1 >\eta$,
at least one of  the numbers
$$
r_2 -|r_1 -1 |,\, r_2 -|r_1 +1 |
$$
must have absolute value $\gtrsim \eta$ and so 
\begin{multline*}
\int_{T_1}\int_{T_2}\frac {1}{\sqrt{\big|\big(r_2 -|r_1 -1 |\big)\big(r_2 -|r_1 +1 |\big)\big|}}\,dr_1 \, dr_2 \lesssim \\
C(\eta ) 
\int_{T_1}\int_{T_2}\Big(\frac {1}{\sqrt{\big| r_2 -|r_1 -1 |\big|}}
+
\frac {1}{\sqrt{\big| r_2 -|r_1 +1 |\big|}}\Big)
\,dr_2 \, dr_1 .
\end{multline*}
Now \eqref{9} follows since the inner integral is $\lesssim |T_2 |^{1/2}$.
$\ \ \ \ \ \ \ \ \ \ \ \ \ \ \ \ \ \ \ \ \ \ \ \ \ \ \ \ \ \ \ \ \ \ \ \ \ \ \ \ \ \ \ \ \ \ \ \ \ \ \ \ \ \ \ \ \ \ \ \ \ \ \ \ \ \ \ \ \ \ \ \ \ \ \ \ \ \ \ \ \ \ \ \ \ \ \ \ \ \ \square$

We return to the project of bounding $|E|$ from below. If $x_1 ,\dots ,x_N$ are in $A$, then 
we estimate, using \eqref{4.75}, \eqref{4.5}, and \eqref{6}, that
\begin{equation}\label{a}
|E|\ge |\cup_j E_{x_j}|\ge NB\mu -C B^{3/2}\sum_{1\le j<k\le N}\frac{1}{|x_j -x_k |^{1/2}}.
\end{equation}
We need to choose the $x_j \in A$ in order to control the sum
\begin{equation*}
\sum_{1\le j<k\le N}\frac{1}{|x_j -x_k |^{1/2}}
\end{equation*}
and we will do this by choosing the $x_j$ independently from the probability space 
$$
\Big(A,\frac{\lambda}{\lambda (A)}\Big).
$$
We reason
\begin{multline*}
\mathbb{E}\Big(\frac{1\ }{|x_j -x_k |^{1/2}}\Big)=
\frac{1}{\lambda (A)^2}\int_A \int_A \frac{1\ }{|x_j -x_k |^{1/2}}\,d\lambda (x_j )\,d\lambda (x_k )\le \\
\frac{1}{\lambda (A)^2}
\Big(\int_A \int_A 1\,d\lambda (x_j )\,d\lambda (x_k )\Big)^{1-1/2\alpha}
\Big(\int_A \int_A \frac{1 }{|x_j -x_k |^{\alpha}}\,
d\lambda (x_j )\,d\lambda (x_k )\Big)^{1/2\alpha}\le 
\frac{C}{ \lambda (A)^{1/\alpha}},
\end{multline*}
where we have used the hypotheses $\alpha >1/2$ and 
\begin{equation*}
\int_{\bbR^2}\int_{\bbR^2} \frac{1\ }{|x_j -x_k |^\alpha} \,d\lambda (x_j )\,d\lambda (x_k )<\infty .
\end{equation*}
Therefore 
\begin{equation*}
\mathbb{E}\Big(\sum_{1\le j<k\le N}\frac{1\ }{|x_j -x_k |^{1/2}}\Big)\lesssim
\frac{N^2}{ \lambda (A)^{1/\alpha}}
\end{equation*}
and so there are $x_1 ,\dots ,x_N \in A$ for which \eqref{a} gives 
\begin{equation*}
|E|\ge |\cup_j E_{x_j}|\ge NB\mu -\frac{C N^2 B^{3/2}}{\lambda (A)^{1/\alpha}}.
\end{equation*}
An appropriate choice of $N$ then yields 
\begin{equation}\label{b}
\mu^2 B^{1/2}\lambda (A)^{1/\alpha }\lesssim |E|,
\end{equation}
which is equivalent to \eqref{4}.
(The $N$ here will be $\gtrsim 1$ so long as $$\mu B^{-1/2}\lambda ^{1/\alpha}\gtrsim 1.$$
If this fails it is easy to see that \eqref{b} holds anyway.)

The proofs of (b) and (c) follow the same general strategy. 
For (b) we will start with \eqref{a} but use the following lemma to choose the points $x_j$.
\begin{lemma}\label{lemma2.5} Suppose $0<\alpha <\alpha ' <\gamma \le d$ and suppose that $\lambda$ is a nonnegative 
compactly supported Borel measure on $\bbR^d$ satisfying the estimate 
\begin{equation}\label{20}
\lambda \big(B(x,\delta )\big)\le C\,\delta^{\alpha '}\ (x\in\bbR^d, \delta >0).
\end{equation}
Then it is possible to choose $x_1 ,\dots ,x_N \in A$ such that 
\begin{equation}\label{21}
\sum_{1\le i<j\le N}\frac{1}{|x_i -x_j |^\gamma}\lesssim 
\lambda (A)^{-\gamma /\alpha}N^{1+\gamma /\alpha},
\end{equation}
where the implied constant depends on $\alpha ,\alpha ' , \gamma, C$, and the size of the support of $\lambda$.

\end{lemma}

{\it Proof of Lemma \ref{lemma2.5}:}
For $k=2, \dots N$, set 
\begin{equation*}
\eta_k =c\Big(\frac{\lambda (A)}{k}\Big)^{1/\alpha }
\end{equation*}
where $c$ is a small positive constant which will depend only on $\alpha$, $\alpha '$, the size of the support of $\lambda$, and $C$ from \eqref{20}: if
$x_1 ,\dots x_{N-1}\in A$, $1\le k\le N-1$, and if we define
\begin{equation*}
A(x_1 ,\dots ,x_k )=\{x\in A: |x-x_j |\ge \eta_j ,\, 1\le j \le k\},
\end{equation*}
then $c$ should be chosen small enough to guarantee that 
\begin{equation}\label{22}
\lambda \big(A(x_1 ,\dots ,x_k )\big)\ge \lambda (A)/2 .
\end{equation}

Now define 
\begin{equation*}
\tilde A_N =\{(x_1 ,\dots ,x_N )\in A^N: x_1 \in A,\, x_2 \in A(x_1 ) ,\dots ,\, x_N \in A(x_1 ,\dots ,x_{N-1})\}
\end{equation*}
and define a probability measure $\tilde \lambda$ on $\tilde A_N $ by
\begin{multline*}
\langle f ,\tilde\lambda \rangle =
\frac{1}{\lambda (A)}\int_A \frac{1}{\lambda \big(A(x_1 )\big)}\int_{A(x_1)}\cdots  \\
\frac{1}{\lambda \big(A(x_1 ,\dots x_{N-1} )\big)}\int_{A(x_1 ,\dots ,x_{N-1})} 
f(x_1 ,\dots x_N )\, d\lambda (x_N )\cdots d\lambda (x_2 )\,d\lambda (x_1 ).
\end{multline*}
If $j<k$, taking an expectation with respect to $\tilde\lambda$ gives 
\begin{multline*}
\bbE \Big(\frac{1}{|x_j -x_k |^\gamma }\Big)\lesssim 
\eta_j^{\alpha -\gamma} \lambda (A)^{-2}\int_{A\times A}\frac{1}{|x_j -x_k |^\alpha }  d(\lambda \times \lambda )(x_k ,x_j )
\lesssim \\
\eta_j^{\alpha -\gamma}\lambda (A)^{-1}\lesssim \lambda (A)^{-\gamma /\alpha}j ^{\gamma /\alpha -1},
\end{multline*}
where we have used \eqref{20} and $\alpha' >\alpha$. Thus 
\begin{equation*}
\bbE \Big(\sum_{1\le j<k\le N}\frac{1}{|x_j -x_k |^\gamma }\Big)\lesssim \lambda (A)^{-\gamma /\alpha}
\sum_{1\le j<k\le N}
j^{\gamma /\alpha -1} \lesssim \lambda (A)^{-\gamma /\alpha}N^{1+\gamma /\alpha}.
\end{equation*}

$\ \ \ \ \ \ \ \ \ \ \ \ \ \ \ \ \ \ \ \ \ \ \ \ \ \ \ \ \ \ \ \ \ \ \ \ \ \ \ \ \ \ \ \ \ \ \ \ \ \ \ \ \ \ \ \ \ \ \ \ \ \ \ \ \ \ \ \ \ \ \ \ \ \ \ \ \ \ \ \ \ \ \ \ \ \ \ \ \ \ \square$

We return to the proof of (b). Taking $\gamma =1/2$, Lemma \ref{2.5} shows that we can choose $x_1 ,\dots ,x_N \in A$ with 
\begin{equation*}
\sum_{1\le j<k\le N}\frac{1}{|x_j -x_k |^{1/2}}\lesssim \lambda (A)^{-1 /2\alpha}N^{1+1 /2\alpha}.
\end{equation*}
This time  \eqref{a} gives 
\begin{equation*}
|E|\ge |\cup_j E_{x_j}|\ge NB\mu -\frac{CN^{1+1 /2\alpha}B^{3/2}}{\lambda (A)^{1/2\alpha}}
\end{equation*}
and an optimal choice of $N$ gives

\begin{equation}\label{c}
\mu^{1+2\alpha} B^{1-\alpha}
\lambda (A)\lesssim |E| .
\end{equation}
(Again, $N$ here will be $\gtrsim 1$ so long as $\mu B^{-1/2}\lambda (A)^{1/2\alpha}\gtrsim 1$,
and if this fails it is easy to see that \eqref{c} holds anyway.) This completes the proof of (b).

The proof of (c) is again very similar. In this case the analogue of \eqref{6} is the estimate 
\begin{equation}\label{24}
|E_{x_1}\cap E_{x_2}|\lesssim \frac{B^2}{|x_1 -x_2 |},
\end{equation}
which we will establish below.
With \eqref{24} and the choice $\gamma =1$ in Lemma \ref{lemma2.5}, \eqref{4.75} leads to
\begin{equation*}
|E|\ge |\cup_j E_{x_j}|\ge NB\mu -\frac{C N^{1+1/\alpha} B^{2}}{\lambda (A)^{1/\alpha}}.
\end{equation*}
One then proceeds as in the proof of (b). We omit the details here and will conclude the proof of Theorem \ref{est3} by sketching a proof of \eqref{24}.
Our starting point for this is the well known (see, e.g., Lemma 1 in \cite{O2}) estimate 
\begin{equation}\label{a2}
|A_\delta (x_1, r_1 )\cap A_\delta (x_2 ,r_2 )|\lesssim \frac{\delta^2}{\delta +|x_1 -x_2|+|r_1 - r_2 |},
\end{equation}
valid for annuli in $\bbR^d$ so long as $d\ge 3$. Assume that $\delta$ is small and that the subsets 
$S(x_i )\subset [r_0 ,R_0 ]$ can be written as disjoint unions of $2\delta$-length intervals:
\begin{equation*}
S(x_i ) =\cup_{j=1}^J [r^i_j -\delta ,r^i_j +\delta ] .
\end{equation*}
Then since $|S(x_i )|\approx B$ we have $J\delta \approx B$. Now
\begin{equation*}
E_{x_i}\subset \cup_{j=1}^J A_\delta (x_i, r^i_{j} )
\end{equation*}
so, by \eqref{a2},
\begin{equation*}
|E_{x_1}\cap E_{x_2}|\le \sum_{j_1 ,j_2} |A_\delta (x_1, r^1_{j_1} )\cap A_\delta (x_2 ,r^2_{j_2} )|
\lesssim \\
 \sum_{j_1 ,j_2} \frac{\delta^2}{\delta +|x_1 -x_2|+|r^1_{j_1} - r^2_{j_2} |}.
 \end{equation*}
For fixed $j_1$ a rearrangement argument and $J\delta\approx B$ show that 
 \begin{equation*}
 \sum_{ j_2} \frac{\delta}{\delta +|x_1 -x_2|+|r^1_{j_1} - r^2_{j_2} |} \lesssim 
 \int_{|x_1 -x_2 |}^{|x_1 -x_2 |+B}\frac{1}{r}\, dr.
\end{equation*}
Thus, using $J\delta\approx B$ again, 
\begin{equation*}
|E_{x_1}\cap E_{x_2}|\lesssim B\log (1+B/|x_1 -x_2 |),
\end{equation*}
more than enough to establish \eqref{24}.

$\ \ \ \ \ \ \ \ \ \ \ \ \ \ \ \ \ \ \ \ \ \ \ \ \ \ \ \ \ \ \ \ \ \ \ \ \ \ \ \ \ \ \ \ \ \ \ \ \ \ \ \ \ \ \ \ \ \ \ \ \ \ \ \ \ \ \ \ \ \ \ \ \ \ \ \ \ \ \ \ \ \ \ \ \ \ \ \ \ \ \square$

\ \ \ 

{\it Proof of Theorem \ref{est4}:} The proof is analogous to the proof of Theorem \ref{est2}. We begin with (a). Fix 
an arbitrary $\eta >0$.
We will show that
\begin{equation*}
\dim \{x\in\bbR^2 : \dim D_x (E) <(2\beta -1)/3 -\eta \}=0.
\end{equation*}
Fix $\alpha$ and $\alpha '$ with $0<\alpha <\alpha' <1/2$. Suppose that $\lambda$ is compactly supported and satisfies \eqref{2.6}. It is enough to show that if $\tilde F\subset\bbR^2$ is compact and $\lambda (\tilde F )>0$, then there is 
some compact $F$ with $F\subset \tilde F$, $\lambda (F)>0$, and $\dim D_x (E)> (2\beta -1)/3 -\eta$ for $\lambda$-a.a.$\,x\in F$. 

Fix a small $\epsilon>0$. Set $\beta '=\beta -\epsilon$. Suppose $p$, $q$, and $s$ 
are defined by \eqref{2.65} with $t=1-\epsilon$, then put 
$$
\rho =\beta ' +\frac{1}{p}(2-\beta ')-\epsilon ,
$$
and finally define $\tau$ by 
$$
\rho -1 =\tau +\frac{1}{s}(1-\tau ).
$$
Observe that if 
$\epsilon$ were $0$ then we would have $\tau =(2\beta -1)/3$. Thus if 
$\epsilon$ is small enough, we do have $\tau >(2\beta -1)/3 -\eta$.
Since $\dim E=\beta$ there is a probability measure $\nu$ on $E$ satisfying 
$$
\nu \big(B(x,\delta )\big)\lesssim \delta^{\beta '} \ (x\in \bbR^2,\, \delta>0).
$$
Now let $\tilde F$ be as above and then choose $r_0$, $R_0$, $r(1)$, $r(2)$,  and $F$ as in the proof of Theorem \ref{est2}. 
It follows from the discussion before Lemma \ref{besov} and the definition of $\rho$ that 
$\nu \ast K_\rho \in L^p (\bbR^ 2)$. It then follows, as in the proof of Theorem \ref{est2}
but using (b) of Theorem \ref{est3} instead of Theorem \ref{est1}, that 
$\dim D_x (E) \ge \tau$ for $\lambda$-a.a.$\,x\in F$.  

The proof of (b) of Theorem \ref{est4} uses (c) of Theorem \ref{est3} instead of (b) of that theorem and is otherwise analogous to the argument above.

$\ \ \ \ \ \ \ \ \ \ \ \ \ \ \ \ \ \ \ \ \ \ \ \ \ \ \ \ \ \ \ \ \ \ \ \ \ \ \ \ \ \ \ \ \ \ \ \ \ \ \ \ \ \ \ \ \ \ \ \ \ \ \ \ \ \ \ \ \ \ \ \ \ \ \ \ \ \ \ \ \ \ \ \ \ \ \ \ \ \ \square$

\begin{equation*}
\end{equation*}

\end{document}